\documentclass{amsart}
  \usepackage{url}
  \usepackage{amssymb}
  
  \input xypic
  \xyoption{all}
  
  \usepackage{amsmath, amsthm,  amsfonts}
 
\let\MathTimeFonts= N 
\let\FormalScriptFont= N

\if Y\MathTimeFonts% 
    \usepackage{mathtime}  
\SetMathAlphabet{\mathit}{normal}{\encodingdefault}{\rmdefault}{m}{it}% 
\fi 
 
\if Y\FormalScriptFont% 
   \usepackage{calrsfs} 
\fi 
 
\makeatletter 
 
\newtheorem{The}{Theorem}[section]

\newtheorem{Rem}[The]{Remark} 
\newtheorem{Lem}[The]{Lemma}

\renewcommand{\theenumi}{\alph{enumi}}

\renewcommand{\p@enumi}{} 
\renewcommand{\p@enumii}{(\theenumi)}

\def\endproof{\relax\ifmmode\expandafter\endproofmath\else 
\unskip\nobreak\hfil\penalty50\hskip.75em\hbox{}\nobreak\hfil\bull 
{\parfillskip= 0pt \finalhyphendemerits= 0 \bigbreak}\fi} 
\def\endproofmath$${\eqno\bull$$\bigbreak} 
\def\bull{\vbox{\hrule\hbox{\vrule\kern3pt\vbox{\kern6pt}\kern3pt\vrule} 
\hrule}}  
 
\def\cancel#1#2{\ooalign{$\hfil#1\mkern1mu/\hfil$\crcr$#1#2$}} 
\def\Dirac{\mathpalette\cancel D}

\newcommand{\ba}{\begin{eqnarray}} 
\newcommand{\na}{\end{eqnarray}} 
\newcommand{\ban}{\begin{eqnarray*}} 
\newcommand{\nan}{\end{eqnarray*}}

\newcommand{\C}{\mathbb{C}} 
 
\newcommand{\R}{\mathbb{R}} 
\newcommand{\Z}{\mathbb{Z}}

\renewcommand{\P}{\mathbb{P}} 
\newcommand{\A}{\mathbb{A}} 
 \newcommand{\V}{\mathbb{V}}
 
 \renewcommand{\cH}{{\mathcal  H}}
\newcommand{\cA}{{\mathcal  A} }
\newcommand{\G}{{\mathcal G}}

\newcommand{\F}{{\mathcal F}}

\newcommand{\U}{{\mathcal U}}
\newcommand{\D}{{\mathcal D}}

\renewcommand{\d}{\partial} 
 
\newcommand{\disp}{\displaystyle}

\title{On the relationship of gerbes to the odd families index theorem} 

 \author[A.L. Carey]{Alan L. Carey}
\address[Alan L. Carey]
  {Mathematical Sciences Institute\\
  Australian National University\\
  Canberra ACT 0200 \\
  Australia}
  \email{acarey@maths.anu.edu.au}

  \author[Bai-Ling Wang]{Bai-Ling Wang}
  \address[Bai-Ling Wang]
  {Mathematical Sciences Institute\\
  Australian National University\\
  Canberra ACT 0200 \\
  Australia}
  \email{wangb@maths.anu.edu.au}
\thanks{This work is  supported by  the Australian
  Research Council Discovery Project DP0449470.}
  
\begin{document}

\begin{abstract}  
The goal of this paper is to 
apply the universal gerbe of \cite{CMi1} and \cite{CMi2}
to give an alternative, simple and more unified view 
of the relationship between index theory
and gerbes. We discuss determinant bundle gerbes \cite{CMMi1}
and the index gerbe of \cite{L} for the case of families
 of Dirac operators on odd dimensional closed manifolds. 
The method also works
 for a family of Dirac operators on odd dimensional 
  manifolds with boundary, for a pair of Melrose-Piazza's 
 $Cl(1)$-spectral sections for a family of Dirac operators on 
 even dimensional closed manifolds with vanishing index in $K$-theory
and, in a simple case, for manifolds with corners.
The common
 feature of these bundle gerbes is that there exists a canonical bundle gerbe
connection whose curving is given by the degree 2 part of the even
eta-form (up to a locally defined exact form) arising from the local family index theorem.
\end{abstract}
\maketitle 
%\Subjclass{57J52;55R65;19K56;58J28.}

%\Keywords{Bundle gerbe; Dixmier-Douady class; local family index theorems; eta forms.}

\tableofcontents

\section{Introduction} 

To the authors' knowledge the subject started with an unpublished
manuscript of Graeme Segal \cite{Seg} who described a substitute,
for families of Dirac operators on odd dimensional manifolds,
for Quillen's determinant line bundle.
Subsequently in \cite{CMu1} Segal's construction was recognised as
defining a bundle gerbe in the sense of  Murray 
\cite{Mur}.
Murray's bundle gerbes are
 differential geometric objects which offer an alternative to
 Brylinski's
description \cite{Bry} of Giraud's gerbes.
In \cite{CMMi1} an explicit formula for the
Dixmier-Douady class of the `determinant bundle gerbe' of \cite{CMu1}
was derived from the odd local index theorem for particular families of
Dirac operators.  We became aware of the existence of a more general
point of view in discussions with Mathai on his work with Melrose \cite{MaM}.
This theory was realised in \cite{L}, where Lott gave a construction
of the higher degree analogs of the determinant line bundle:
the `indexgerbe' and defined connective structures on them, 
for the case of a family
of Dirac operators on odd-dimensional closed spin manifolds. He
also recognised the role of eta forms.

The so-called
universal gerbe was introduced in \cite{CMi1}\cite{CMi2} where 
Carey and Mickelsson 
analysed the obstruction to obtaining a second quantization for a smooth 
family of Dirac
operators on an odd dimensional spin manifold
with boundary. Explicit computations of the
Dixmier-Douady class for the universal gerbe were also given.
With hindsight \cite{L} may be seen as an example of
the universal gerbe for the case of families of Dirac operators.
 The index gerbe is studied from a slightly different point of view by
by Bunke in \cite{Bu1}.

In this paper we use the universal gerbe to simplify the 
discussion of \cite{L} and to provide a way to extend \cite{L}.
We also  show that the universal gerbe can be used 
to construct other  examples of geometrically
interesting gerbes such as that for  
a family of Dirac operators on 
odd dimensional manifolds with boundary and for a pair of Melrose-Piazza 
$Cl(1)$-spectral sections for a family of Dirac operators on 
even dimensional closed manifolds with vanishing index in $K$-theory. 
These latter
gerbes exhaust all bundle gerbes corresponding to the Dixmier-Douady classes
in the image of the Chern character map on the $K^1$-group of the underlying
manifold. 
In all these examples of gerbes, there exists a canonical bundle gerbe
connection whose curving, up to an exact 2-form, 
is given by the degree 2 part of the even 
eta-form arising from the local family index theorem.
This enables us to propose an approach to family 
index theory manifolds with corners
by using the bundle gerbe whose curving is given by 
the degree 2 part of the eta-form.
 
For a single manifold with corners, Fredholm perturbations of
Dirac type operators and their index are thoroughly studied by
Loya and Melrose in \cite{LM}. A more general and ambitious
approach to these latter questions is contained in \cite{Bu2}
which motivated us to develop the point of view of this paper.

As this paper is largely about gaining a unifying perspective
we need to review earlier work.
Section 2 summarises what is needed from the theory of bundle gerbes.
Section 3 describes briefly
the  family index theorems developed by 
Bismut\cite{Bis2}, Bismut-Freed\cite{BF}, 
Bismut-Cheeger\cite{BC2} and Melrose-Piazza
\cite{MP1}\cite{MP2}. This Section is required to set up notation.
In Section 4, we recall the construction of the 
determinant bundle gerbe in \cite{CMMi1}
and study its geometry in Theorem \ref{DD:gauge}.
In Section 5, we give a complete proof of the existence theorem
(Theorem \ref{universal:gerbe}) of a canonical
bundle gerbe class over a base manifold associated to any element in 
$K^1$ of the base manifold using the universal gerbe.
The results in this Section make explicit 
some folklore of the subject.

 In Section 6, we apply the 
 constructions in Section 5 
to the index gerbe associated to a family of
 Dirac operators on an odd dimensional manifold with or without boundary. 
 This gives a simple construction of the gerbe discussed in 
 \cite{L} and \cite{Bu1}. 
 Two new examples of bundle gerbes are obtained, one for
 a family of Dirac operators on odd dimensional 
  manifolds with boundary,  and the other for a pair of Melrose-Piazza
 $Cl(1)$-spectral sections for a family of Dirac operators on 
 even dimensional closed manifolds with vanishing index in $K$-theory. 
 Our main observation in this paper is that the universal
 gerbe suggests an approach to understanding
 the local family index theory (in particular the transgression of eta forms)
for manifolds with codimension two corners.
Our approach, if it can be fully developed, may provide an alternative
to  the treatment of index theory on manifolds with boundary and corners
in \cite{Bu2} and \cite{BuP}.

\noindent{\bf Acknowledgement} ALC 
acknowledges advice from P. Piazza on spectral sections
and also thanks V. Mathai for conversations on his joint work
with R. Melrose related to the theme analysed in this paper.  BLW wishes
to thank Xiaonan Ma for many invaluable comments and his explanations 
of eta forms.
\section{Review of bundle gerbes}

Our reference for this section is \cite{Mur}.
Take a smooth surjective submersion $\pi: Y\to M$ and let $Y^{[2]}$ 
be the fiber product associated with $\pi$
with the obvious groupoid structure. 
A bundle gerbe $\G$ over $M$, as  defined in  \cite{Mur},
consists of $\pi:  Y\to M$, and a principal $U(1)$-bundle $\G$
over $Y^{[2]}= Y\times_\pi Y$ together with a groupoid multiplication
on $\G$, which is compatible with the natural  groupoid multiplication
on $Y^{[2]}$. The multiplication is represented by an isomorphism
   \ba
   \label{bg:prod}
   m:\qquad  \pi^*_1\G \otimes \pi^*_3\G   \to \pi^*_2\G
   \na
   of principal $U(1)$-bundles over $Y^{[3]} = Y\times_\pi Y\times_\pi Y$, 
where $\pi_i$, i=1,2, 3,  are the three natural projections from
   $Y^{[3]}$ to $Y^{[2]}$ 
obtained by omitting the entry in position $i$ for $\pi_i$.
   
In \cite{MurSte} it is shown that the natural equivalence relation on bundle gerbes is that of stable equivalence. Within each stable equivalence 
class we may choose a local representative
by letting $\{U_\alpha\}$ be a good cover.
Then
the submersion $\pi$ admits a local section $s_\alpha$ over $U_\alpha$.
The local bundle gerbe can be described by a family of local
$U(1)$-bundles $\{\G_{\alpha\beta}\}$ over $\{U_{\alpha\beta}= U_\alpha\cap U_\beta\}$,
the pull-back of $\G$ via $(s_\alpha, s_\beta)$.
In this picture the bundle gerbe multiplication (\ref{bg:prod})
is given by an isomorphism  
\ba\label{cech}
\phi_{\alpha\beta\gamma}: \qquad \G_{\alpha\beta}\otimes \G_{\beta\gamma} \cong \G_{\alpha\gamma},
\na
over $U_{\alpha\beta\gamma}$, such that $ \phi_{\alpha\beta\gamma}$
is associative as a groupoid multiplication.
A related local picture of gerbes 
was introduced by Hitchen \cite{Hit}.

A C\v{e}ch cocycle $\{f_{\alpha\beta\gamma}\}$ can be obtained from the
isomorphism (\ref{cech}) by choosing a section $s_{\alpha\beta}$
of $\G_{\alpha\beta}$, i.e.,
\[
\phi_{\alpha\beta\gamma}(s_{\alpha\beta}\otimes s_{\beta\gamma}) = 
f_{\alpha\beta\gamma}\cdot s_{\alpha\gamma}
\]
for a $U(1)$-valued function $f_{\alpha\beta\gamma}$ over 
$U_{\alpha\beta\gamma}$. The equivalence class
of $\{f_{\alpha\beta\gamma}\}$ doesn't depend on the choice of
local sections $\{s_{\alpha\beta}\}$,
and represents the
Dixmier-Douady class of the bundle gerbe in
\[
H^2_{\text{C\v{e}ch}}(M, \underline{U(1)}) \cong H^3(M, \Z),
\]
where $\underline{U(1)}$ is the sheaf of continuous $U(1)$-valued functions over $M$. 

The geometry of bundle gerbes, bundle gerbe connections and their curvings
is studied in \cite{Mur}. On the corresponding local bundle gerbe 
 a bundle gerbe connection is a family of
$U(1)$-connections $\{\nabla_{\alpha\beta}\}$ on $\{\G_{\alpha\beta}\}$
 which is compatible
with the isomorphism (\ref{cech}), i.e., $A_{\alpha\beta} = \nabla_{\alpha\beta}
 s_{\alpha\beta}
/s_{\alpha\beta}$ satisfies 
\ba\label{compatible}
A_{\alpha\beta} + A_{\beta\gamma} +A_{\gamma\alpha} = 
f_{\alpha\beta\gamma}^{-1}d f_{\alpha\beta\gamma}.
\na
A curving then is a locally defined 2-form $B_\alpha$
 such that, 
\ba\label{curving}
B_\alpha- B_\beta  = dA_{\alpha\beta},
\na
over $U_\alpha\cap U_\beta$. We remark that 
in string theory applications this system $\{B_\alpha\}_\alpha$ of local 
two-forms is called the $B$-field. We will use this terminology here.
The $B$-field
 is unique,  up to locally defined exact 2-forms. If we choose
another representative for the $B$-field
  $\{B_\alpha + d C_\alpha\}_\alpha$, then we need to modify
 the gerbe connection by adding $\{C_\alpha-C_\beta\}$ to 
$\{A_{\alpha\beta}\}$.
This has no effect on
 the compatibility condition (\ref{compatible}).  This local description of
connection and   curving for a
local bundle gerbe defines an element 
\[
[(f_{\alpha\beta\gamma}, A_{\alpha\beta},
B_\alpha)]
\]
in the degree three Deligne cohomology group
$H^3_{Del}(M, \D^3)$.  The  curvature of $B$-field is given by 
$\{dB_\alpha\}$ and we note that 
\[
dB_\alpha = dB_\beta, \qquad  \text{on $U_{\alpha\beta}$}
\]
is a globally defined 3-form which  represents the image of the
Dixmier-Douady class in $H^3(M, \R)$ up to a factor of $\disp{\frac{i}{2\pi}}$.
This 3-form is the gerbe analogue of the Chern class 
of line bundles with connections. In this paper, we often suppress the
normalisation factor
and identify the bundle gerbe curvature with the differential form representing
the image of the Dixmier-Douady class in $H^3(M, \R)$.

In fact, though we do not use it here,
the  $B$-field of a bundle gerbe $\G$ over $M$
is a Hopkins-Singer differential 
cocycle as defined in \cite{HopSin} whose equivalence class is  the  corresponding 
degree three  Deligne cohomology class. 
One has in fact the following commutative diagram
 \[
 \xymatrix{
 H^3_{Del}(M, \D^3) \ar[d] \ar[r]& H^3(M, \Z)\ar[d]\\
 \Omega^3_\Z(M) \ar[r] & H^3(M, \R)}
 \]
 where $\Omega^3_\Z(M)$ is the space of closed 3-forms on $M$ with
 integral periods. 
We also note that there have been other applications 
of bundle gerbes to various aspects string theory
and quantum field theory:
\cite{BoMa}\cite{BCMMS} \cite{CMu1}\cite{CMMi1}\cite{CMMi2}
\cite{CMi1}\cite{CMi2} \cite{EM}\cite{GR}
\cite{Mic}.

\noindent{\bf Important remark}.
There is a subtlety here relating to the   gauge transformations
of  the $B$-field. To illustrate this subtlety, 
given a globally defined closed 2-form $\omega$
(not necessarily of integer periods), with respect to a cover $\{U_\alpha\}$ of
$M$ we can write  $\omega |_{U_\alpha} = d\theta_\alpha$ as an exact 2-form by
the Poincare Lemma. Then  on the one hand,
one can see that 
\[
(1, 0, \omega |_{U_\alpha})
\]
represents a trivial degree three
Deligne class in $H^3_{Del}(M, \D^3)$ if and only if the closed 2-form
$\omega$ has integer periods. This claim follows from the exact sequence
\[
0\longrightarrow \Omega^2(M)/\Omega^2_\Z(M) \longrightarrow
H_{Del}^3(M, \D^3) \longrightarrow H^3(M, \Z) \longrightarrow 0.
\]
 In this sense,  the gauge transformation of the
$B$-field is given by line bundles with connection on $M$. On the other hand, 
\[
(1, \theta_\alpha -\theta_\beta, \omega |_{U_\alpha}) =
 (1, \theta_\alpha -\theta_\beta, d\theta_\alpha)
\]
 always represents
a trivial degree three Deligne class in $H^3_{Del}(M, \Z)$, as 
$(1, \theta_\alpha -\theta_\beta, d\theta_\alpha)$ is a coboundary element.
 In this paper,
a curving on a bundle gerbe is always  defined up to a locally
defined exact 2-form in this latter sense so that the corresponding
 $B$-field is well-defined up to a degree three Deligne coboundary
 term of form $(1, \theta_\alpha -\theta_\beta, d\theta_\alpha)$. Here
 $d\theta_\alpha$ may not be a globally defined exact 2-form. 
 We hope that this remark will clarify any confusion when we talk
 about curvings on a bundle gerbe, such
  as the degree 2 part of the even eta-form up to an
 exact 2-form. 
 
\section{Review of local family index theory}

In this Section, we briefly review local index theory and eta forms
and set up notation.
Local family index theorems refine the cohomological
family index theorem of Atiyah-Singer, 
while eta forms transgress the local family index theorem
at the level of differential forms.  The main reference for this section
is Bismut's survey paper \cite{Bis1} and the references therein.

Let $\pi: X \to B$ be a smooth fibration over a closed smooth manifold
$B$, whose fibers are diffeomorphic to a 
compact, oriented, spin manifold $M$. Let $(V, h^V, \nabla^V)$ be 
an Hermitian vector
bundle over $X$ equipped with a unitary connection.

Let $g^{X/B}$ be a metric on the relative tangent bundle $T(X/B)$ 
(the vertical tangent subbundle of $TX$), which
is fiberwise a product 
metric near the boundary if the boundary of $M$ is non-empty. In the latter
case, $\partial X \to B$ is also a fibration with compact fiber
$\partial M$.
Let $S_{X/B}$ 
be the corresponding spinor bundle associated to the spin structure on
$(T(X/B), g^{X/B})$.

 Let $T^HX$ be a smooth vector subbundle of $TX$, (a horizontal vector
 subbundle of $TX$),  such that 
 \ba\label{horizontal}
 TX = T^HX \oplus T(X/B),
 \na
 and such that if $\partial M$ is non-empty,  $T^HX|_{\d X} \subset T(\d X)$.
 Then $T^H(\d X) = T^HX|_{\d X}$ is a horizontal subbundle of $T(\d X)$.

 Denote by $P^v$ the projection of $TX$ onto the vertical tangent bundle
 under the decomposition (\ref{horizontal}).   
Theorem 1.9 of \cite{Bis2} shows that $(T^HX, g^{X/B)}$ determines a 
canonical Euclidean
 connection $\nabla^{X/B}$ as follows. Choose a metric $g^X$ on $TX$ such that
 $T^HX$ is orthogonal to $T(X/B)$ and such that $g^{X/B}$ is 
the restriction of 
 $g^X$ to $T(X/B)$. Then 
 \[
 \nabla^{X/B} = P^v \circ \nabla^X
 \]
 where $\nabla^X$ is the Levi-Civita connection on $(TX, g^X)$.

 Write 
 \[
 \hat A (x) = \disp{\frac {x/2}{sinh(x/2)}},\ \ 
 Ch (x) = exp(x).
 \]
 Then the  characteristic classes can be represented
 by closed differential forms on $X$ as
 \[
 \hat A (T(X/B), \nabla^{X/B}) = det^{1/2} \bigl(\hat A ( \disp{\frac{ i(\nabla^{X/B})^2}
 {2\pi}}) \bigr),
 \]
 \[
 Ch(V, \nabla^V) =Tr\bigl( Ch(\disp{\frac{ i(\nabla^{V})^2}{2\pi}})\bigr).
 \]

For any $b\in B$, let $\Dirac_b$ be the Dirac operator acting on 
$C^\infty (X_b, (S_{X/B}\otimes V)|_{X_b})$, where
$X_b= \pi^{-1}(b)$ is the fiber of $\pi$ over $b$. 
Then $\{\Dirac_b\}_{b\in B}$ 
is a  family of 
elliptic  operators on the fibers of  $\pi$, acting fiberwisely on an infinite dimensional
vector space 
$\C^\infty
(X_b, (S_{X/B}\otimes V)|_{X_b}),
$
parametrized by $b\in B$.
With respect to the inner product from
the Hermitian metrics on the spinor bundle $S_{X/B}$ and $V$,  $\{\Dirac_b\}_{b\in B}$ 
can be viewed as a family of unbounded operators  
$$\Dirac_b:  L^2(X_b, (S_{X/B}\otimes V)|_{X_b}) \longrightarrow 
L^2(X_b, (S_{X/B}\otimes V)|_{X_b}).$$

\begin{Rem}\begin{enumerate}
\item If  the fiber  of $\pi: X\to B$ is even dimensional,
then the spinor bundle $S_{X/B}= S_{X/B}^+ \oplus S_{X/B}^-$ is
 $\Z_2$-graded, and 
 \[
\Dirac_b = \left( \begin{array}{ccc}
0 & \Dirac_{+, b}^*\\
\Dirac_{+, b}&0
\end{array} \right).
\]
  Then $\Dirac_{+, b} (1+  \Dirac_{+, b}^*\Dirac_{+, b})^{-\frac 12}$
is a continuous family of  Fredholm operators,   
which  defines a $K^0$-element, denoted by  $Ind (\Dirac_+)$ in $K^0(B)$. 
The Atiyah-Singer index theorem becomes
\ba\label{AS}
Ch\bigl(Ind (\Dirac_+)\bigr) = \pi_*[\hat A (T(X/B), \nabla^{X/B}) Ch(V, 
\nabla^V)] \in H^{even}(B, \R)
\na
where $\pi_*$ is given by the integration along the fibers
on the level of differential forms.

 \item If  the fiber  of $\pi: X\to B$ is odd dimensional,
 the family $\{\Dirac_b\}_{b \in B}$ defines a  $K^1$-element
$ Ind (\Dirac) \in K^1(B)$. 
  The Atiyah-Singer index theorem becomes
 \ba\label{AS:odd}
 Ch(Ind (\Dirac )) =  \pi_*[\hat A(T(X/B), \nabla^{X/B}) Ch(V,\nabla^V)] 
\in H^{odd}(B, \R).
\na
\end{enumerate}
\end{Rem}

Now assume that  the fibers  of the fibration $\pi: X \to B$ are diffeomorphic
to an even dimensional manifold 
$M$ with non-empty boundary $\d M$, then the family of boundary
 Dirac operators $\{\Dirac^{\d M}_b\}_{b\in B}$
  for the boundary fibration $\d\pi: \d X \to B$ defines a zero class 
  in $K^1(B)$ by cobordism invariance of the Atiyah-Singer index theory. Hence, there
  exists  an even eta-form $\tilde {\eta}_{even}$ on $B$ satisfying the transgression formula
 \[
  d \tilde {\eta}_{even} =  (\d\pi)_*\bigl(
\hat A (T(\d X/B), \nabla^{\d X/B}) Ch (V, \nabla^V)\bigr),
\]
 $\tilde {\eta}_{even}$ is unique up to an exact odd form on $B$.
 
 Melrose-Piazza introduced the notion of a
 spectral section for
the family of boundary  Dirac operators $\{\Dirac^{\d M}_b\}_{b\in B}$ in \cite{MP1},
where a spectral
section is a smooth family of self-adjoint pseudo-differential projections
\[
P_b: \qquad L^2 (\d X_b , (S_{\d X/B}\otimes V)|_{\d X_b})
\to L^2 (\d X_b , (S_{\d X/B}\otimes V)|_{\d X_b}),
\] such that there is 
$R> 0$ for which  $\Dirac^{\d M}_b s = \lambda s $ implies that
\[
\begin{array}{l}
P_b s = s \qquad \qquad \lambda > R,
\\
P_b s =0 \qquad \qquad \lambda< -R.
\end{array}\]
For a fixed spectral section $P$, then Melrose-Piazza  constructed 
 an even eta form $\tilde{\eta}_{even}^{P}$
depending only on 
$\{\Dirac^{\d M}_b\}_{b\in B}$ and the spectral section $P$.  
With the boundary condition given by a spectral section $P$,
they also established the local
family index theorem for  the family of Dirac operators 
$\{\Dirac^M_{+, P; b}\}_{b\in B}$ in 
  (\cite{MP1}).

\begin{The} \label{MP:1} Let $\pi: X \to B$ be a smooth fibration with fiber diffeomorphic
to an even dimensional spin manifold $M$ with boundary.
The family of Dirac operators 
$\{\Dirac^M_{+, P; b}\}_{b\in B}$ has a well-defined index
$Ind (\Dirac^M_{+, P}) \in K^0(B)$, with the Chern character  given by
\[
Ch \bigl(Ind (\Dirac_{+, P}^M)\bigr) = \disp{\int_M} 
\hat A (T(X/B), \nabla^{X/B}) Ch (V, \nabla^V) -\tilde{\eta}_{even}^P \in
H^{even}(B, \R).
\]
\end{The}

When the fibers of the fibration $\pi: X \to B$ are diffeomorphic to
an odd dimensional manifold $M$
with non-empty boundary, the family of boundary Dirac operator
$\{\Dirac^{\d M}_b\}_{b\in B}$ consists of self-adjoint, elliptic, 
odd differential operators acting on 
\ba\label{polarized}
  L^2 (X_b, (S^+_{\d X/B}\otimes V)|_{\d X_b}) \oplus
L^2 (X_b, (S^-_{\d X/B}\otimes V)|_{\ X_b}).
\na

In order to get a continuous family of self-adjoint Fredholm operators, 
 Melrose and Piazza introduced a $Cl(1)$-spectral section
$P$ and constructed an odd eta form $\tilde {\eta}^P_{odd}$ which, modulo exact forms,
depends only on $\{\Dirac^{\d M}_b\}_{b\in B}$ and the $Cl(1)$-spectral section 
(\cite{MP2}). The $Cl(1)$-condition (where $Cl(1)$ denotes the Clifford algebra $Cl(\R)$)
is given by
\[
c(du) \circ P  + P \circ c(du) = c(du),
\]
with $u$ being the inward normal coordinate near the
boundary $\d M$. 

\begin{The}\label{MP:2}
(\cite{MP2}) Let $\pi: X \to B$ be a smooth fibration with fiber diffeomorphic
to an odd dimensional spin manifold $M$ with boundary. 
With the boundary condition given by a $Cl(1)$-spectral section $P$,  
the family of Dirac operators 
$\{\Dirac^M_{P; b}\}_{b\in B}$ has a well-defined index
\[
Ind (\Dirac^M_{P}) \in K^1(B),
\] with the Chern character form given by
\ba
Ch \bigl(Ind (\Dirac_{ P}^M)\bigr) = \disp{\int_M} 
\hat A (T(X/B), \nabla^{X/B}) Ch (V, \nabla^V) -\tilde{\eta}_{odd}^P 
\in H^{odd}(B, \R).
\na
\end{The}

\section{Determinant bundle gerbe}

In \cite{CMMi1}, Carey-Mickelson-Murray used the determinant bundle gerbe
to study
the bundle of fermionic Fock spaces parametrized by
connections on principal $G$-bundles on odd dimensional
manifolds (for $G$ a compact Lie group). 
We briefly review this construction here using the local family index theorem.

Let $M$ be a closed, oriented, spin manifold 
with a Riemannian metric $g^M$. Let $S$ be the
spinor bundle over $M$ and let $(V, h^V)$ be a Hermitian vector bundle.
Denote by $\cA$ the space of unitary connections on $V$ and by
$\G$ the based gauge transformation group, 
that is, those gauge transformations
fixing the fiber over some fixed point in $M$. With proper regularity 
on connections,
the quotient space $\cA/\G$ is a smooth Frechet manifold. 

Let $\lambda \in \R$, and let
\[
U_\lambda =\{A\in \cA| \lambda \notin spec(\Dirac_A)\},
\]
where the Dirac operator $\Dirac_A$ acts on $C^\infty (M, S\otimes V)$. Let 
$H$ be the space of square integrable sections of $S\otimes V$. For any
$A\in U_\lambda$, there is a uniform polarization:
\[
H= H_+(A, \lambda) \oplus H_-(A, \lambda)
\]
given by  the spectral decomposition of $H$ with respect to $\Dirac_A-\lambda$ into the eigenspaces corresponding to positive and negative spectrum.
 Denote by $P_{H_\pm(A, \lambda})$ the 
orthogonal projection onto $H_\pm (A, \lambda)$.

Fix $\lambda_0 \in \R$ and  a reference connection
$A_0\in \cA$ such that the  Dirac operator $\Dirac_{A_0}-\lambda_0$ is
invertible and denote by $P_{H_\pm (A_0, \lambda_0)}$ the 
orthogonal projection onto $H_\pm (A_0, \lambda_0)$.
Over $U_\lambda$, there exists a complex line bundle
$Det_\lambda$, which is essentially the determinant line bundle for the 
even
dimensional Dirac operators on $[0, 1]\times M$ coupled to some path of connections
in $\cA$ with respect to the
generalized Atiyah-Singer-Patodi boundary condition (\cite{APS}).
 The path of connections $A(t)$  
can chosen to be  any smooth path connecting $A_0$ and 
$A\in U_\lambda$, for simplicity, we take 
\ba\label{A(t)}
\A= A(t) = f(t)A + (1-f(t)) A_0
\na
for $t \in [0, 1]$ where $f$ is zero on $[0,1/4]$, equal to one on
$[3/4, 1]$ and interpolates smoothly between these values on $[1/4,3/4]$.
The generalized  Atiyah-Singer-Patodi boundary condition is determined by
the orthogonal spectrum projection 
\[
P_\lambda = P_{H_+(A_0, \lambda_0)} \oplus P_{H_-(A, \lambda)},
\]
that is,  the even dimensional Dirac operator on $[0, 1]\times M$ acts on 
the spinor fields whose boundary components belong to
$H_-(A_0, \lambda_0)$ at $\{0\}\times M$ and $H_+(A, \lambda)$ at
$\{1\}\times M$.

Note that $P_\lambda$ can be thought as a spectral section in the sense of
Melrose-Piazza \cite{MP1} for the family of Dirac operators
\[
\{\Dirac_{A_0} \oplus \Dirac_A\}_{A\in U_\lambda},
\]
on $\{0\}\times M\sqcup \{1\}\times M$.
 Then the eta form $\tilde{\eta}^{P_\lambda}_{even}$ associated with
 $\{\Dirac_{A_0} \oplus \Dirac_A\}_{A\in U_\lambda}$ and the
 spectral section $P_\lambda$ is an even differential form on
 $U_\lambda$, which is unique up to an exact form.

 Denote by $\Dirac_\A^{P_\lambda}$ the even dimensional Dirac operator
$\Dirac_\A$ with respect to the APS-type boundary condition $P_\lambda$.
Then the even dimensional Dirac operator
\[
\Dirac_\A^{P_\lambda}=
\left( \begin{array}{ccc}
0 & \Dirac_{\A,-}^{P_\lambda}\\
\Dirac^{P_\lambda}_{\A,+}&0
\end{array} \right)
\]
 is a Fredholm operator, whose
determinant line is given by 
\[
Det_\lambda (A) =  \Lambda^{top}(Ker \Dirac_{\A, +}^{P_\lambda})^* \otimes 
\Lambda^{top}(Coker \Dirac_{\A, +}^{P_\lambda}).
\]

The family of Fredholm operators $\{\Dirac_\A^{P_\lambda}\}$, parametrized
by $A\in U_\lambda$ defines a determinant line bundle 
  over $U_\lambda$, which is given by
\[
Det_\lambda = \bigcup_{A\in U_\lambda} Det_\lambda(A).
\]
The determinant line bundle
$Det_\lambda$ can be equipped with a Quillen metric and a Bismut-Freed
unitary connection whose curvature can be calculated by
the local family index theorem.

  Over $U_{\lambda\lambda'} = U_\lambda \cap U_{\lambda'}$, there exists a 
  complex line   bundle $Det_{\lambda\lambda'}$ such that
  \[
  Det_{\lambda\lambda'} = Det_\lambda^* \otimes Det_{\lambda'}.
  \]
These local line bundles $\{Det_{\lambda\lambda'}\}$ 
over $\cA$, form a bundle gerbe as established 
in \cite{CMMi1}. As $\cA$ is contractible this bundle gerbe is trivial
so the interest in \cite{CMMi1} is lies in what happens 
when one takes the induced bundle gerbe
on $\cA/\G$. 

In the following, we will apply the local family index theorem
to study the geometry of this determinant bundle gerbe.
To apply the local family index theorem, we should restrict ourselves
 to a smooth finite dimensional
submanifold of $\cA$. For convenience however, we will formally work on 
 the infinite dimensional manifold $\cA$  directly. 
 
 Consider the trivial fibration $[0, 1] \times M \times \cA$
  over $\cA$ with fiber  $[0, 1] \times M$ an even dimensional manifold 
  with boundary. 
 Over $[0, 1] \times M \times \cA$, there is a Hermitian
 vector bundle $\V$  which is the pull-back bundle of $V$. 
 
 There is a universal
 unitary connection on $\V$, also denoted by $\A$, whose
 vector potential at $(t, x, A)$ is given by 
 \[
 \A (t, x) = A(t)(x),
 \]
 where $A(t)$ is given by (\ref{A(t)}). Denote by $Ch(\V, \A)$
 the Chern character of $(\V, \A)$.

  Now we can state the following theorem regarding the geometry of
  the bundle gerbe   over $\cA/\G$ constructed in \cite{CMMi1}.
  
  \begin{The}\label{DD:gauge}
  The local line bundles $\{Det_{\lambda\lambda'}\}$ descend to 
local line bundles
  over $\cA/\G$, which in turn define a local bundle gerbe
  over $\cA/\G$. Moreover, the induced unitary connection and the even eta 
  form $\tilde{\eta}_{even}^{P_\lambda}$ 
  (up to an exact 2-form) descend to a bundle gerbe connection and
  curving on the local bundle gerbe over $\cA/\G$, whose bundle gerbe curvature is 
  given by the differential form 
  \[
  \disp{
  \bigl( \disp{\int_M} \hat A(TM, \nabla^{TM})Ch(\V, \A)\bigr)_{(3)}}.
  \]
    \end{The}
    \begin{proof}
  Equip $Det_\lambda$ with the Quillen metric and its
unitary connection, then
 its first Chern class is represented by the degree 2 part of the
 following  differential form:
 \[
 \disp{\int_{[0, 1]\times M}} \hat A(TM, \nabla^{TM})Ch(\V, \A)
 -\tilde{\eta}_{even}^{P_\lambda}.
 \]
In this formula
$\nabla^{TM}$ is the Levi-Civita connection on $(TM, g^M)$,
$\hat A(TM, \nabla^{TM})$ represents the $\hat A$-genus of $M$, and
  $\tilde{\eta}_{even}^{P_\lambda}$ is the even eta form on $U_\lambda$
associated
  to the family of boundary Dirac operators and the spectral section
  $P_\lambda$. Note that $\tilde{\eta}_{even}^{P_\lambda}$ (modulo exact forms)
  is uniquely determined, see Theorem \ref{MP:1}.
  
  The induced connection on $Det_{\lambda\lambda'}$ implies
  that its first Chern class is given by
  \[
  \bigl( \tilde{\eta}_{even}^{P_\lambda} - \tilde{\eta}_{even}^{P_{\lambda'}}
  \bigr)_{(2)}.
  \]
  Note that the eta form $\tilde{\eta}_{even}^{P_\lambda}$ is unique up to
  an exact form (in a way analogous to the $B$-field), hence, $\bigl( \tilde{\eta}_{even}^{P_\lambda} - \tilde{\eta}_{even}^{P_{\lambda'}}
  \bigr)_{(2)}$ is a well-defined element in $H^2(M, \R)$.
   From the transgression formula for the eta forms and Stokes formula, 
  we know that  over $U_{\lambda\lambda'}$,
  \[
  d (\tilde{\eta}_{even}^{P_\lambda})_{(2)} -
  d (\tilde{\eta}_{even}^{P_{\lambda'}})_{(2)}
  = \bigl( \disp{\int_M} \hat A(TM, \nabla^{TM})Ch(\V, \A)\bigr)_{(3)}.
  \]
  Here use the fact  that the contribution from $\{0\} \times M$ vanishes
  as a differential form on $\cA$, as we fixed a connection $A_0$ 
  over $\{0\} \times M$. We also use the same notation
  $(\V, \A)$ to denote the Hermitian vector bundle 
  over $M \times \cA$ and the universal unitary connection $\A$.
  
  As the gauge group acts on $\cA$ and $Det_{\lambda}$ covariantly, by 
  quotienting out $\G$, we obtain the bundle gerbe over $\cA/\G$ described 
  in the theorem.
  \end{proof}
  
  \begin{Rem}  The Dixmier-Douady class in Theorem \ref{DD:gauge}
is often non-trivial.
 For example note that when $dim M < 4$, there is no
  contribution from the $\hat A$-genus
and in  particular, for $dim M =1$ or $3$, using
  the Chern-Simons forms, non-triviality was proved by  explicit
calculation  in \cite{CMMi1}.
  \end{Rem}
  
The above construction can be generalized to the fibration case as in the local
family index theorem for odd dimensional manifolds with boundary of Section 3.
When  the fibers are closed odd dimensional spin manifolds, 
this leads to the index gerbe as discussed by Lott (\cite{L}). 
In the next section, we discuss the universal bundle gerbe as 
in \cite{CMi1}\cite{CMi2}, which provides a unifying viewpoint for those
bundle gerbes constructed from various determinant
line bundles.

\section{The universal bundle gerbe}

Let $H$ be an infinite dimensional separable complex Hilbert space. Let
$\F^{a.s}_*$ be the space of all self-adjoint Fredholm operators on $H$ with
positive and negative essential spectrum. With the norm topology on
$\F^{a.s}_*$,  Atiyah-Singer \cite{AS1} showed that
$\F^{a.s}_*$ is a representing space for the $K^1$-group, that is, for any 
closed manifold $B$,
\[
K^1(B) \cong [ B, \F^{a.s}_*]
\]
the homotopy classes of continuous maps from $B$ to $\F^{a.s}_*$. 
As $\F^{a.s}_*$ is homotopy equivalent to $\U^{(1)}$, the group of unitary
automorphisms $g: H \to H$ such that $g-1$ is trace class, we obtain
\[
K^1(B) \cong [B, \U^{(1)}].
\]

In \cite{CMi1}, a universal bundle gerbe was constructed on $\U^{(1)}$ with the
Dixmier-Douady class given by the basic 3-form
\ba
\label{basic:form}
\disp{\frac{1}{24\pi^2}}Tr (g^{-1}dg)^3,
\na
the generator of $H^3(\U^{(1)}, \Z)$. 

For any compact Lie group, this basic
3-form gives the so-called  basic gerbe. 
We prefer however to call it the universal bundle
gerbe in the sense that many examples of  bundle gerbes on 
smooth manifolds are obtained by pulling back this universal bundle gerbe
via  certain smooth maps from $B$ to $\U^{(1)}$.

Note that the odd Chern character of $K^1(B) = [B, \U^{(1)}]$ is given by
\ba\label{odd:chern}
Ch([g]) = \sum_{n\geq 0} \disp{\left[(-1)^n\frac{n!}{(2\pi i)^{n+1} (2n+1)!}
Tr\bigl( (g^{-1}dg)^{2n+1}\bigr) \right]},
\na
for a smooth map $g: B \to \U^{(1)}$ representing a $K^1$-element
$[g]$ in $K^1(B)$. This odd Chern character formula was proved in \cite{Getz}.

The following theorem was established in \cite{CMi1}
as the first obstruction to obtaining a 
second quantization for a smooth family of Dirac operators (parametrized by
$B$) on an odd dimensional spin manifold. 
By a second quantization,  we mean an irreducible representation of the
canonical anticommutative relations (CAR) algebra, the
complex Clifford algebra $Cl(H\oplus \bar{H})$, which is compatible
with the action of the quantized Dirac operator. Such a 
representation is given
by the Fock space associated to a polarization on $H$. For
a bundle of Hilbert spaces over $B$, a continuous polarization always
exists locally, but not necessarily globally. 
This leads to a bundle gerbe over $B$.

\begin{The}
\label{universal:gerbe}
 For any $K^1$-element $[g] \in K^1(B)$ represented by a smooth map 
 $g: B \to \U^{(1)}$. There exists a canonical construction of a 
bundle gerbe over $B$ with connection and curving whose
 bundle gerbe curvature is given by the degree
3 part of 
\[
\disp{\frac{1}{24\pi^2}}Tr (g^{-1}dg)^3.
\]
\end{The}
\begin{proof}
The bundle gerbe we are after is essentially the pull-back
bundle gerbe from the universal bundle gerbe over
$\U^{(1)}$ under the smooth map $g: B \to \U^{(1)}$.

Consider the polarized Hilbert space $\cH = L^2(S^1, H)$ with the polarization
$\epsilon$ given by the Hardy decomposition.
\[
\cH = L^2(S^1, H) = \cH_+ \oplus \cH_-.
\]
 That is we take
the polarisation given by splitting into positive and negative Fourier
modes. Then the smooth based loop group
$\Omega \U^{(1)}$ acts naturally on $\cH$. From \cite{PS}, we see that 
\ba\label{inclusion}
\Omega \U^{(1)} \subset \U_{res}(\cH, \epsilon).
\na
Here $\U_{res}(\cH, \epsilon)$ is the restricted unitary group of $H$ with respect
to the polarization $\epsilon$, those $g\in U(\cH)$ such that
the off-diagonal block of $g$ is Hilbert-Schmidt.
It was shown in \cite{CMi1} that the inclusion 
$\Omega \U^{(1)} \subset \U_{res}(\cH, \epsilon)$ is a homotopy equivalence.

We know that the holonomy map from the space of 
connections on a trivial $\U^{(1)}$-bundle
over $S^1$ provides a model for the universal $\Omega \U^{(1)}$-bundle. Hence,
$\U^{(1)}$ is a classifying space for $\Omega \U^{(1)}$. 

   From $K^1(B) \cong [B, \U^{(1)}]$, we conclude that elements
in $K^1(B)$ are in  one-to-one correspondence with 
isomorphism classes of principal  $\U_{res}(\cH, \epsilon)$-bundles over $B$. 

Associated to 
the basic three form (\ref{basic:form}) on $\U^{(1)}$, is the universal
gerbe realized 
as the lifting bundle gerbe \cite{Mur} associated to the central extension
\ba
\label{extension}
1\to U(1) \to \hat{\U}_{res} \to \U_{res}(\cH, \epsilon) \to 1.
\na
We recall some of the theory behind this fact.
First note that the fundamental representation of $ \hat{\U}_{res}$
acts on the Fock space (see \cite{PS})
\[
\F_{\cH} = \Lambda (\cH_+) \otimes \Lambda(\bar{\cH}_-),
\]
associated to the polarized Hilbert space $\cH= \cH_ \oplus \cH_-$. 
This gives rise to a homomorphism $\U_{res}(\cH, \epsilon) \to \P U (\F_\cH)$
which induces a $\P U (\F_\cH)$ principal bundle from the 
$\U_{res}(\cH, \epsilon)$ principal bundle
associated to the $K^1$-element $[g]$. 
Denote by ${\mathcal P}_{g}$
the resulting  $\P U (\F_\cH)$ principal bundle corresponding to a smooth map
$g: B \to \U^{(1)}$ representing $[g]$.
The corresponding lifting bundle gerbe (as defined in \cite{Mur})
is the canonical bundle gerbe over $B$ for which we are looking.

Recall that the lifting bundle gerbe associated 
to ${\mathcal P}_{g}$
can be described locally as follows
\cite{Mur}. Take a local trivialization of ${\mathcal P}_g$
with respect to a good cover of $B=\bigcup_{\alpha} U_\alpha$. Assume over
$U_\alpha$, the trivialization is given by a local section 
\[
s_\alpha: \qquad U_\alpha \to {\mathcal P}_g|_{U_\alpha},
\]
such that the transition function over $U_{\alpha\beta} = U_\alpha \cap U_\beta$
is given by a smooth function 
\[
\gamma_{\alpha\beta}:\qquad U_{\alpha\beta} \to  \U_{res}(\cH, \epsilon)
\subset \P U (\F_\cH).
\]
We can define a local Hermitian line bundle 
${\mathcal L}_{\alpha\beta}$ over $U_{\alpha\beta}$ as 
follows. First pull back
the principal $U(1)$ bundle (\ref{extension}) using $\gamma_{\alpha\beta}$.
Then construct the 
associated line bundles $\{ {\mathcal L}_{\alpha\beta}\}$
over the double intersections $U_{\alpha\beta}$.
This family of  local Hermitian line bundles 
$\{ {\mathcal L}_{\alpha\beta}\}$ defines a bundle gerbe over $B$
with multiplication obtained from the multiplication in $\hat{\U}_{res}$.

In \cite{CMi1} an Hermitian connection on 
$\{ {\mathcal L}_{\alpha\beta}\}$ is given which defines
a bundle gerbe connection, together  with a  choice of a  curving such that 
the bundle gerbe curvature is given by
\[ \disp{\frac{1}{24\pi^2}}Tr (g^{-1}dg)^3.
\]
\end{proof}

For an Hermitian vector bundle $(V, h^V)$ over an oriented, closed, 
spin manifold $M$ with a Riemannian metric
denote by $\cA$ the space of unitary connections on $V$ and by
$\G$ the based gauge transformation group as in section 4.

In \cite{CMi2}, an explicit smooth map  $\tilde g: \cA/\G \to \U^{(1)}$ was constructed:
firstly assign to any $A\in \cA$ a unitary operator in $\U^{(1)}$ by
\[
A\to g(A)=-exp(i\pi \disp{\frac{\Dirac_A}{|\Dirac_A|+ \chi (|\Dirac_A|)}}),
\]
where $\chi$ is any positive smooth exponentially decay function on $[0, \infty)$
with $\chi(0)=1$. As $g(A)$ is not gauge invariant,
\[
g(A^u) = u^{-1}g(A) u,
\]
for any $u\in \G$,  it doesn't define a smooth map 
from $\cA/\G$ to $ \U^{(1)}$. In order
to get a gauge invariant map, a global section 
for the associated $U(H)$-bundle
$\cA\times_\G U(H)$ is needed.
This exists due to the contractibility of $U(H)$. This
section is given by a smooth map $r: \cA \to U(H)$ such that 
$r(A^u) = u^{-1}r(A)$. Then the required map $\tilde g: \cA/\G \to \U^{(1)}$ 
is given by $\tilde{g}(A) = r(A)^{-1}g(A )r(A).$

Let $B$ be a smooth submanifold of $\cA/\G$. The restriction of
$\tilde g$ to $B$ defines an element in $K^1(B)$, which is exactly the 
family of Dirac operators over $B$ associated to the universal connection
$\A$ on $\V$ (see section 3 for the definition). 

  From the local family index theorem, we know that  
the curvatures satisfy  
  \ba\label{curvatures}
  \disp{
  \bigl( \disp{\int_M} \hat A(TM, \nabla^TM)Ch(\V, \A)\bigr)_{(3)}} 
  =  \disp{\frac{1}{24\pi^2}}Tr (\tilde g^{-1}d\tilde g)^3.
 \na
Note that the left hand side was established in \cite{CMMi1} and the right hand
side in \cite{CMi2}.
This relation fixes the image of
 the Dixmier-Douady class  in $H^3(B, \R))$. It follows that
the bundle gerbe over $B\subset \cA/\G$ in Theorem \ref{universal:gerbe}
 is stably isomorphic to the determinant
bundle gerbe constructed in Theorem \ref{DD:gauge}, up to taking a product
with 
a bundle gerbe with torsion Dixmier-Douady class. 
(The notion of products for bundle gerbes is covered in the original
paper of \cite{Mur}). This is because the Dixmier-Douady class of a product 
is the sum of the Dixmier-Douady classes and a bundle gerbe with torsion 
Dixmier-Douady class can be equipped with a connection and curving whose
bundle gerbe curvature is trivial. Such a  bundle gerbe is often 
called  flat.

\begin{Rem} There is in fact an equivalent but different picture
for the universal bundle gerbe.
If a  $K^1$-element in $K^1(B)$ is represented by a smooth
family of self-adjoint Fredholm operators $\{T_b\}_{b\in B}\in\F^{a.s}_*$
on a Hilbert space $H$, there is a 
canonical principal  $\P U (\F_H)$-bundle over $B$ constructed as follows. 
First we work over $F^{a.s}_*$
          by covering it with open sets of form
\[
U_\lambda = \{  T\in \F^{a.s}_* | \lambda \in spec (T)\},
\]
where $\lambda\in \R$. Over $U_\lambda$, there is a
polarization 
\[
H = H^-_{T,\lambda} \oplus H^+_{T,\lambda},
\]
varying continuously in $T$
given by the spectral decomposition of $H$ into eigenspaces 
of $T\in U_\lambda$ corresponding
to eigenvalues greater or less than $\lambda$. Denote by $P_{T,\lambda}$ the
orthogonal projection onto $H^-_{T,\lambda}$. 
For $\lambda\neq \lambda'$, $P_{T,\lambda}$
and $P_{T,\lambda'}$ are related by 
conjugation by an element $g_{\lambda\lambda'}(T)$ 
 of $\hat{\U}_{res}$ from the lifting of a well-defined element 
\[
\tilde g_{\lambda\lambda'}(T)\in \U_{res},
\]  
where the copy of  
$\hat{\U}_{res}$ we are using is specified by a reference spectral projection
of $T$. Now there is a Fock representation corresponding to 
this reference spectral projection and hence,
using the inclusion  $\U_{res} \subset \P U (\F_H)$
a well-defined element 
\[
\tilde g_{\lambda\lambda'}(T)\in\P U (\F_H).
\]
The family of functions $\{\tilde g_{\lambda\lambda'}\}$
defined on overlaps $U_{\lambda\lambda'}$
 can be used as transition functions to construct
a principal  $\P U (\F_H)$-bundle over $F^{a.s}_*$.
We get a corresponding bundle over $B$ using the pullback construction.

The above argument  also 
gives a universal determinant bundle gerbe $\mathcal D$ over $F^{a.s}_*$
by using the construction in Section 4 via determinant line bundles
$Det_{\lambda\lambda'}$
over intersections $U_{\lambda\lambda'}$.
Recall from \cite{CMMi1} that for $\lambda<\lambda'$,
the line over $T$ is
 $Det_{\lambda\lambda'}(T)$:
the highest exterior power of the vector space spanned
by eigenvectors of $T$ corresponding to eigenvalues between
$\lambda$ and $\lambda'$.
The Dixmier-Douady class of this determinant bundle gerbe
can be determined exactly.
To see this we observe that the Dixmier-Douady
class of the pull back of  $\mathcal D$
to $\U_1$ using an explicit homotopy equivalence from \cite{AS1} must be a multiple,
say $n$, of the fundamental class on $\U_1$.
To determine which multiple just take
 $B$ to be $S^3$
and choose a  family $\{T_b\}_{b\in S^3}\in\F^{a.s}_*$
to represent the generator of $H^3(\F^{a.s}_*, \Z)\cong \Z$
(an explicit example is given in the last section of \cite{CMMi1})
Under the explicit homotopy equivalence from \cite{AS1}, we know
that the resulting  continuous map $S^3 \to \U_1$ also defines the
generator of $H^3(\U_1 , \Z)\cong \Z$.   The 
fact that $n=1$ follows from (\ref{curvatures})
as in the case of $S^3$ the curvature suffices to determine the
Dixmier-Douady class.
\end{Rem}

\section{Index gerbe as induced from the universal bundle gerbe}

In this Section, we will construct the index gerbe associated to a family 
of Dirac operators on an odd dimensional manifold with or without boundary.
To do this we need to recall some
additional standard material using the notation of Section 3.

Let   $\pi: X \to B$ be a smooth fibration 
with even dimensional fibres.
With a choice of spectral section $P$ for 
the family of boundary  Dirac operators $\{\Dirac^{\d M}_b\}_{b\in B}$, 
we have a family of Fredholm operators $\{\Dirac^M_{+, P; b}\}_{b\in B}$
over $B$.  There is a 
determinant line bundle, denoted by $Det(\Dirac^M_{+, P})$,
over $B$ given by 
\[
Det(\Dirac^M_{+, P})_b = \Lambda^{top}(Ker \Dirac^M_{+,P; b})^* \otimes \Lambda^{top}
(Coker \Dirac^M_{+,P; b}).
\]
 Using zeta determinant regularization, 
as in  \cite{Qu}\cite{BF}\cite{Woj}\cite{SW}, 
a Hermitian metric and a unitary connection $\nabla^{\Dirac^M_{+, P}}$ 
can be constructed on $Det(\Dirac^M_{+, P})$ such that
\ba\label{c1:det}
c_1(Det(\Dirac^M_{+, P}), \nabla^{\Dirac^M_{+, P}}) 
= \bigl[\pi_*\bigl(\hat A (T(X/B), \nabla^{X/B}) Ch (V, \nabla^V)\bigr) -
\tilde{\eta}_{even}^P\bigr]_{(2)},
\na
Note that the eta form $\tilde{\eta}_{even}^P$ is defined modulo exact forms, the above
equality (\ref{c1:det}) holds modulo exact 2-forms. 

Now we let $\pi: X\to B$ have
 closed even-dimensional fibers partitioned into 
two codimension zero submanifolds, 
$M= M_0\cup M_1$, joined 
 along a codimension 1 submanifold $\d M_0 = -\d M_1$. 
Assume that
the metric $g^{X/B}$ is of product type near the collar neighborhood of the
separating submanifold. Let $P$ be a spectral section 
for $\{\Dirac^{\d M_0}_b\}_{b\in B}$, then $I-P$ is a spectral section for
$\{\Dirac^{\d M_1}_b\}_{b\in B}$. 
Scott showed (\cite{Sco}) that the determinant
line bundles for these three families of 
$\{\Dirac^M_{+; b}\}_{b\in B}$, $\{\Dirac^{M_0}_{+, P; b}\}_{b\in B}$
and $\{\Dirac^{M_1}_{+, I-P; b}\}_{b\in B}$ satisfy the following gluing formulae
as Hermitian line bundles:
\ba\label{split:det}
Det(\Dirac^M_{+}) \cong Det (\Dirac^{M_0}_{+, P}) \otimes 
Det (\Dirac^{M_1}_{+, I-P}),
\na
moreover the splitting formulae for the curvature of $\nabla^{\Dirac^M_{+}}$
implies that 
\ba\label{split:curv}
c_1(Det(\Dirac^M_{+}), \nabla^{\Dirac^M_{+}}) = 
c_1(Det(\Dirac^{M_0}_{+, P}), \nabla^{\Dirac^{M_0}_{+, P}})
+ c_1(Det(\Dirac^{M_1}_{+, I-P}), \nabla^{\Dirac^{M_1}_{+, I-P}}).
\na
With these facts in hand we move on to the main results of this paper.

\subsection{Index gerbe from a family of Dirac operators on a closed  
manifold.}

Let $\pi: X \to B$ be a smooth fibration over a closed smooth manifold
$B$, whose fibers are diffeomorphic to a 
compact, oriented, odd dimensional spin manifold $M$.
 
Let $g^{X/B}$ be a metric on the relative tangent bundle $T(X/B)$ and
let $S_{X/B}$ be the  spinor bundle associated to 
$(T(X/B), g^{X/B})$. 
 Let $T^HX$ be a horizontal vector subbundle of $TX$. Then
 $(T^HX, g^{X/B})$ determines a connection $\nabla^{X/B}$ on
 $T(X/B)$  as in section 3.
Let $(V, h^V, \nabla^V)$ be a Hermitian vector
bundle over $X$ equipped with a unitary connection.

The family of Dirac operators
 $\{\Dirac_b\}_{b \in B}$ defines a $K^1$-element
 \[
  Ind (\Dirac )\in K^1(B),
  \]
  with $Ch(Ind (\Dirac ) )
= \pi_*\bigl(\hat A (T(X/B), \nabla^{X/B}) Ch (V, \nabla^V)\bigr)
  \in H^{odd}(B, \R)$
  as given by the Atiyah-Singer index theorem (Theorem (\ref{AS:odd})).
  
  Then Theorem
 \ref{universal:gerbe} provides a canonical bundle gerbe
 $ \G^M $  over $B$ with bundle gerbe connection and curving
  whose curvature is given by 
  \[
 \pi_*\bigl(\hat A (T(X/B), \nabla^{X/B}) Ch (V, \nabla^V)\bigr)_{(3)}.
 \]

  Now we can prove the following theorem, which was obtained in  
  \cite{L} by Lott using a different method.
  
  \begin{The}\label{gerbe:closed}
  Let $\pi: X \to B$ be a smooth fibration with fibers   diffeomorphic to 
closed  odd dimensional spin manifolds $M$and
$V$ be a Hermitian vector
bundle over $X$ equipped with a unitary connection  $\nabla^V$.
Then the  associated family of Dirac operators defines a canonical
 bundle gerbe  $\G^M$ over $B$ equipped with a Hermitian metric and a unitary 
 gerbe connection
 whose curving (up to an exact form) is given by the locally defined eta form  
 such that its bundle gerbe curvature is given by 
 \[
 \bigl(\disp{\int_M\hat A (T(X/B), \nabla^{X/B}) Ch (V, \nabla^V)}\bigr)_{(3)}.
 \]
 \end{The}
 \begin{proof}
  Cover  $B$ by $U_\lambda$ (for $\lambda \in \R$)  such that 
  \[
  U_\lambda =\{ b \in B | \Dirac_b -\lambda \text{ is invertible} \}.
  \]
  Over $U_\lambda$, the bundle of Hilbert spaces of square integrable
  sections along the fibers has a continuous polarization:
  \[
  H_b = H^+_b(\lambda) \oplus H^-_b(\lambda),
  \]
  the spectral decomposition with respect to $\Dirac_{P, b}^M-\lambda$ into the
  positive and negative eigenspaces.
  Denote by $P_\lambda$ the continuous family of orthogonal projections
  $\{P_{\lambda, b} \}_{b\in U_\lambda}$ onto $H^+_b(\lambda)$ along the above 
continuous
  polarization over $U_\lambda$.
  
  Note that the restriction of $Ind(\Dirac)$ on $U_\lambda$ is trivial by the
  result of Melrose-Piazza in \cite{MP1},  as $\Dirac$
  over $U_\lambda$ admits a spectral section.  So
  over $U_\lambda$, its Chern character form is exact:
  \ba
  \label{transgression:lambda}
   \pi_*\bigl(\hat A (T(X/B), \nabla^{X/B}) Ch (V, \nabla^V)\bigr) 
   = d \tilde{\eta}_{even}^{P_\lambda}
   \na
   where $\tilde{\eta}_{even}^{P_\lambda}$, unique up to an exact form,
    is the even eta form on $U_\lambda$, 
    associated to $\{\Dirac^M_b\}_{b\in U_\lambda}$ and  the spectral section $P_\lambda$.

  For $\lambda \geq \lambda'$, over $U_{\lambda\lambda'} = U_\lambda \cap U_{\lambda'}$, 
  consider the 
  fibration $[0, 1]\times X \to B$ with even dimensional fibers, the boundary
  fibration has two components $\{0\} \times X$ and $\{1\}\times X$. 
  The spectral section can be chosen to be
  \ba
  \label{spec:section}
  P_{\lambda\lambda'} = P_{\lambda} \oplus (I-P_{\lambda'}).
  \na
  
  Associated with this spectral section $P_{\lambda\lambda'}$, there
  exists an even eta form $\tilde{\eta}_{even}^{P_{\lambda\lambda'}}$
   on $U_{\lambda\lambda'}$, modulo
  exact forms, depending only on $P_{\lambda\lambda'}$ and the family of
  boundary Dirac operators. From the definition of the
  eta form, we have 
  \[
  \tilde{\eta}_{even}^{P_{\lambda\lambda'}}
  = \tilde{\eta}^{P_\lambda}_{even} - 
 \tilde{\eta}_{even}^{P_{\lambda'}},
 \]
 which should be  understood modulo exact even forms on $U_{\alpha\beta}$.
  Then the family of Dirac operators (with boundary condition given by
  the  spectral section   $P_{\lambda\lambda'}$) denoted by 
  \[
  \{\Dirac^{[0, 1] \times M}_{P_{\lambda\lambda'}, b}\}_{b\in U_{\lambda\lambda'} },
  \]
  is a smooth family of Fredholm operators, which
  has a well-defined index $Ind (\Dirac^{[0, 1] \times M}_{P_{\lambda\lambda'}})$
  in $K^0(U_{\lambda\lambda'})$. 
  
  The corresponding determinant line bundle over $U_{\lambda\lambda'}$,
  denoted by $Det_{\lambda\lambda'}$, is a Hermitian line bundle 
  equipped with the Quillen metric and the Bismut-Freed 
   unitary connection. Its first Chern class is given by
  the local family index formula (\ref{c1:det}):
  \ba
  \bigl[\pi_*\bigl(\hat A (T([0, 1]\times X/U_{\lambda\lambda'}), 
  \nabla^{X/B}) Ch (V, \nabla^V)\bigr) -   \tilde{\eta}_{even}^{P_{\lambda\lambda'}}
  \bigr]_{(2)}.
  \na

  Note that in this situation, the contribution
  from  the characteristic  class 
  \[
  \pi_*\bigl(\hat A (T([0, 1]\times X/U_{\lambda\lambda'}), 
  \nabla^{X/B}) Ch (V, \nabla^V)\bigr)
  \]
 vanishes, as it has no component in the $[0, 1]$-direction.  Hence, the first
 Chern  class of $Det_{\lambda\lambda'}$ is represented by the form 
 \[
 \bigl(\tilde{\eta}_{even}^{P_{\lambda\lambda'}}\bigr)_{(2)} =
 \bigl( \tilde{\eta}^{P_{\lambda'}}_{even} - 
 \tilde{\eta}_{even}^{P_{\lambda}}\bigr)_{(2)}.
 \]
 Again these eta forms are well-defined only modulo locally defined
  exact 2-forms.   
 This then implies that the even eta form
 \[
 \bigl(\tilde{\eta}_{even}^{P_\lambda}\bigr))_{(2)},
 \]
 which is only  defined over $U_\lambda$, 
  is the curving (up to an exact 2-form)  
  for  the Bismut-Freed connection on $Det_{\lambda\lambda'}$.
  Hence, the gerbe curvature is uniquely determined and  given by
  \[ d(\tilde{\eta}^{P_\lambda}_{even})_{(2)} - 
  d(\tilde{\eta}_{even}^{P_{\lambda'}})_{(2)}= 
  \bigl(\disp{\int_M\hat A (T(X/B), \nabla^{X/B}) Ch (V, \nabla^V)}\bigr)_{(3)}.
  \]
  \end{proof}

\subsection{Index gerbe from a family of Dirac operators on a  manifold with boundary}

Now we assume that the  fibration $\pi: X \to B$ has fiber diffeomorphic to 
an  odd dimensional $Spin$ manifold $M$
with non-empty boundary.

A $Cl(1)$-spectral section 
$P$ for the family of boundary Dirac operators
$\{\Dirac^{\d M}_b\}_{b\in B}$ provides a well-defined index
for  the family of self-adjoint Fredholm operators $\{\Dirac^M_P\}$:
\[
Ind (\Dirac_{ P}^M) \in K^1(B),
\]
with $Ch(Ind (\Dirac ) )= \pi_*\bigl(\hat A (T(X/B), \nabla^{X/B}) Ch (V, \nabla^V)
-\tilde{\eta}^P_{odd}\bigr)
  \in H^{odd}(B, \R)$  as given by  the local family index theorem in 
Theorem \ref{MP:2}.

 Then Theorem  \ref{universal:gerbe} provides a canonical bundle gerbe 
 $ \G_P^{M}$ over $B$ with a bundle gerbe connection and  curving
 whose curvature  is given by:
 \ba
 \label{DD:MP2}
 \bigl( \pi_*(\hat A (T(X/B), \nabla^{X/B}) Ch (V, \nabla^V) )
 -\tilde{\eta}^P_{odd}\bigr)_{(3)}.
 \na
 Here we should emphasize that $\tilde{\eta}^P_{odd}$ 
is an odd eta form associated
 to a perturbation of $\{\Dirac^{\d M}_b\}_{b\in B}$ by a family of 
self-adjoint 
 smoothing operators $\{A_{P, b}\}_{b\in B}$ such that 
 $\Dirac^{\d M}_b + A_{P, b}$ is invertible as in \cite{MP2}. 
 Note that, as previously,
 we always define  $\tilde{\eta}^P_{odd}$ up to an exact form.

We aim to understand the bundle gerbe connection and its
curving such that the bundle gerbe curvature is given by (\ref{DD:MP2}).
Note that the argument for the previous case can't be applied here, as
now $[0, 1] \times M$ is a manifold with corners near $\d M$, and 
to our knowledge the
even eta form for a family of  Dirac operators on odd dimensional manifolds
with boundary which transgresses the odd local index form (\ref{DD:MP2})
has not been found. 
For  a manifold with corners, Fredholm perturbations of Dirac operators 
and their index formulae have been developed by Loya-Melrose \cite{LM} 
but for a family of Dirac operators on a manifold
with corners up to codimension two their theory does not apply.

Instead, we will apply the theory of bundle gerbes 
to find this even eta form which
transgresses the odd local index form (\ref{DD:MP2}) 
and discuss its implications for local
family index theory for a family of Fredholm operators on a 
manifold with corners of a particular type,
$[0, 1]\times M$. 

 For the fibration $\pi: X \to B$, whose fibers are
  odd dimensional manifolds
with non-empty boundary, the 
Melrose-Piazza $Cl(1)$-spectral section
$P$ for  the family of boundary Dirac operators
$\{\Dirac^{\d M}_b\}_{b\in B}$ gives rise to a
 family of self-adjoint Fredholm operators 
$\{\Dirac^M_{P, b}\}_{b\in B}$ with discrete
 spectrum.
 
 Cover $B$ by $U_\lambda$ with $\lambda\in \R$ such that
 \[
  U_\lambda =\{ b \in B | \Dirac^M_{P,b} -\lambda \text{ is invertible} \}.
  \]

\begin{Lem}  Over $U_\lambda$, the family of self-adjoint Fredholm operators
  $\{\Dirac^M_{P}\}$ has trivial index in $K^1(U_\lambda)$.
  \end{Lem} 
  \begin{proof} This follows from
  the fact that, over $U_\lambda$, we have a uniform polarization of 
the Hilbert space of 
  square integrable sections along the fibers with boundary condition given by
  the $Cl(1)$-spectral section $P$,
  \[
  H_{P, b} = H^+_{P, b}(\lambda) \oplus H^-_{P, b}(\lambda),
  \]
  the spectral decomposition with respect to $\Dirac_{P, b}^M-\lambda$ into the
  positive and negative spectral subspaces. 
Denote by $P_\lambda$ the smooth family of
   orthogonal projections
  $\{P_{\lambda, b} \}_{b\in U_\lambda}$ onto $H^+_{P, b}(\lambda)$ 
  along the above uniform   polarization over $U_\lambda$. Over $U_\lambda$,
  $\Dirac^M_{P}$ admits a spectral section $\{P_{\lambda, b} \}_{b\in U_\lambda}$,
  by the   result of Melrose-Piazza in \cite{MP1}, 
  $\{\Dirac^M_{P}\}$ has trivial index in $K^1(U_\lambda)$. 
  \end{proof} 
  
  Hence, there exists an even form $\eta_{even}^{P, \lambda}$ on $U_\lambda$,
  unique up to an exact form, 
satisfying the following transgression formula over
  $U_\lambda$:
  \ba\label{eta:new}
  d\eta^{P, \lambda}_{even} =\pi_*(\hat A (T(X/B), \nabla^{X/B}) Ch (V, \nabla^V) )
 -\tilde{\eta}^P_{odd}.
 \na
 
 We remark that these forms $\{\eta^{P, \lambda}_{even}\}$ should
in fact be
  even eta forms {\em a la} Melrose-Piazza's
 construction in (\cite{MP2})  associated with the family of self-adjoint
  Fredholm operators   $\{\Dirac^M_{P, b}\}_{b\in U_\lambda}$ and the spectral section
  $P_\lambda$ over $U_\lambda$. We do not have a proof however.
  
  From the abstract bundle gerbe theory in \cite{Mur}, we know that this
  locally defined even form $\bigl(\eta_{even}^{P, \lambda}\bigr)_{(2)}$
   on $U_\lambda$ is the
  curving, up to an exact form on $U_\lambda$, for the bundle gerbe induced
  from the universal gerbe as in Theorem \ref{universal:gerbe}. Hence, 
we have established
  the following Theorem.

  \begin{The} let $\pi: X \to B$ be a fibration  with fiber diffeomorphic to 
an  odd dimensional $Spin$ manifold $M$
with non-empty boundary. Let $P$ be a $Cl(1)$-spectral section 
such that $\{\Dirac^M_P\}$ is 
 a continuous  family of self-adjoint Fredholm operators
 parametrized by $B$, then the associated bundle gerbe $\G^{M}_P$ over $B$
 can be equipped with a bundle   gerbe connection and curving 
 $\{(\eta^{P, \lambda}_{even})_{(2)}\}$
 with the  curvature 
  given by
  \[
  d( \eta^{P,\lambda}_{even})_{(2)} = 
  \bigl(\disp{\int_M\hat A (T(X/B), \nabla^{X/B}) Ch (V, \nabla^V)}-\tilde{\eta}^P_{odd}
  \bigr)_{(3)}.
  \]
  \end{The}

\subsection{Index gerbe from a pair of $Cl(1)$-spectral sections}

  From Proposition 4 in \cite{MP2}, we know that if $P_1$ and 
  $P_2$ are two $Cl(1)$-spectral sections for $\{\Dirac^{\d M}_b\}_{b\in B}$, then
  the Atiyah-Singer suspension operation \cite{AS1} defines a difference
  element
  \[
  [P_2-P_1] \in K^1(B),
  \]
  such that 
  \[
  Ind(\Dirac^M_{P_1})-Ind(\Dirac^M_{P_2}) = [P_2-P_1] \in K^1(B),
  \]
  whose Chern character is given by $\tilde{\eta}^{P_2}_{odd}-\tilde{\eta}^{P_1}_{odd}$.
  It was also shown in Proposition 12 in \cite{MP2} that for a fixed $P_1$, as $P_2$
  ranges over all $Cl(1)$-spectral sections, $[P_2-P_1]$ exhausts $ K^1(B)$.
   
  Apply Theorem  \ref{universal:gerbe} again to $[P_2-P_1]$, we obtain a 
  canonical bundle gerbe 
 $ \G_{P_2P_1} $  associated to $[P_2-P_1] \in K^1(B)$, with a bundle
 gerbe connection and curving whose curvature
  is given by
 \[
 \bigl(\tilde{\eta}^{P_2}_{odd}-\tilde{\eta}^{P_1}_{odd}\bigr)_{(3)}.
 \]
 Moreover, we have the following stable isomorphism relating the bundle gerbes associated to
 families of Dirac operators on odd-dimensional manifolds with boundary 
 equipped with two different $Cl(1)$-spectral sections for the family of boundary
 Dirac operators
 \[
 \G^M_{P_1} \cong \G^M_{P_2} \otimes \G_{P_2P_1}.
 \]

Take a smooth fibration $X$ over $B$, with closed
 odd-dimensional fibers partitioned into two codimensional 0 submanifolds 
$M= M_0\cup M_1$ along a codimension 1 submanifold $\d M_0 = -\d M_1$. Assume that
the metric $g^{X/B}$ is of product type near the collar neighborhood of the
separating submanifold. Let $P$ be a $Cl(1)$-spectral section 
for $\{\Dirac^{\d M_0}_b\}_{b\in B}$, then $I-P$ is a $Cl(1)$-spectral section for
$\{\Dirac^{\d M_1}_b\}_{b\in B}$. Then we have the following splitting formula
for the canonical bundle gerbe $\G^M$ obtained in Theorem \ref{gerbe:closed}:
\[
\G^M \cong \G^{M_0}_P\otimes \G^{M_1}_{I-P}.\]


\small 
\begin{thebibliography}{9999999} 


\bibitem[Ati]{Ati}M.~Atiyah, {\em $K$-theory},
(W.A.~Benjamin, New York, 1967).

 

\bibitem[APS]{APS} M.F. Atiyah, V.K. Patodi, I.M. Singer, {\em Spectral 
       asymmetry and Riemannian geometry}, I,II,III; 
       Math. Proc. Cambridge Phil. Soc. 77 (1975) 43-69; 78 (1975) 405--432; 
       79 (1976) 71-99. 
       
\bibitem[AS1]{AS1}     M.F. Atiyah and I.M. Singer, {\em  Index theory for skew-adjoint 
  Fredholm operators.  }  Inst. Hautes tudes Sci. Publ. Math. No. 37, 1969 5--26.
   
\bibitem[AS2]{AS2}  M.F. Atiyah and I.M. Singer, {\em Dirac operators coupled to
vector
potentials.} Natl. Acad. Sci. (USA) \bf 81, \rm 2597 (1984)

\bibitem[BGV]{BGV} Berline, N., Getzler, E., Vergne, M., {\em Heat kernals and the
Dirac operators.}, 1992.

\bibitem[Bis1]{Bis1} J.M. Bismut, {\em  Local index theory, eta invariants and 
holomorphic torsion: a survey.}  Surveys in differential geometry, Vol. III
 (Cambridge, MA, 1996), 1--76,
Int. Press, Boston, MA, 1998. 

\bibitem[Bis2]{Bis2} Bismut, J.-M. {\em The index Theorem for families of 
Dirac operators: two heat equation proofs.}  Invent. Math.  83, 91--151 (1986).

\bibitem[BC1]{BC1} J.M. Bismut, J. Cheeger, {\em $\eta$-invariants and
     their adiabatic limits}. J. Amer. Math. Soc. 2 (1989), no. 1, 33--70.
     
\bibitem[BC2]{BC2}J. M. Bismut, J. Cheeger,  {\em
Families index for manifolds with boundary,
     superconnections and cones I and II.}
     J. Funct. Anal.  89 (1990) 313--363, and  90 (1990) 306--354.
\bibitem[BF]{BF} Bismut, J.-M., Freed D. {\em The analysis of elliptic families. II. 
Eta invariants and the holonomy theorem.}  Comm. Math. Phys.,  107, 103--163 (1986).


\bibitem[Boss]{Boss}B. BooÐ-Bavnbek, K.  Wojciechowski,  {\em
  Elliptic Boundary Problems for Dirac Operators.} Boston: Birkh"user, 1993 

 \bibitem[BCMMS]{BCMMS} P. Bouwknegt, A. Carey, V. Mathai, M.  Murray, D.
 Stevenson,{\em Twisted $K$-theory and 
 $K$-theory of bundle gerbes.}  Comm. Math. Phys.  228  (2002),  no. 1, 17--45.  
  
  
 \bibitem[BoMa]{BoMa}
P. Bouwknegt and V. Mathai,
{\it $D$-branes, $B$-fields and twisted $K$-theory},
J. High Energy Phys. {\bf 03} (2000) 007, hep-th/0002023.
 
 
 \bibitem[Bry]{Bry} J-L, Brylinski, {\em Loop spaces, characteristic classes and 
  geometric  quantization,} Progress in Mathematics, 107. Birkh"user Boston, 1993. 
 
 \bibitem[Bu1]{Bu1} U. Bunke,  {\em Transgression of the index gerbe.}
   Manuscripta Math.  109  (2002),  no. 3, 263--287.
   
\bibitem[Bu2]{Bu2}
 U. Bunke {\em Index theory, eta forms, and Deligne cohomology}, preprint.

 \bibitem[BuP]{BuP}
U. Bunke and J. Park {\em Determinant bundles, boundaries and surgery}
J. Geom. Phys. {\bf 52} (2004) 28-43. 
\bibitem[CMu1]{CMu1} A.L. Carey and M.K. Murray.:
{\em Faddeev's anomaly and bundle gerbes   }.
Letters in Mathematical Physics, 
{\bf 37} (1996)

\bibitem[CMMi1]{CMMi1}
 Alan Carey, Jouko Mickelsson, Michael Murray
{\em Index Theory, Gerbes, and Hamiltonian Quantization}
 Commun.Math.Phys. 183 (1997) 707-722
 
 \bibitem[CMMi2]{CMMi2}
Alan Carey, Jouko Mickelsson, Michael Murray
{\em Bundle Gerbes Applied to Quantum Field Theory}
 Rev.Math.Phys. {\bf 12} (2000) 65-90
 
 \bibitem[CMi1]{CMi1}
 Alan Carey, Jouko Mickelsson
{\em A gerbe obstruction to quantization of fermions on odd dimensional manifolds}
Lett.Math.Phys. {\bf 51} (2000) 145-160

\bibitem[CMi2]{CMi2} A.L.Carey and J. Mickelsson
{\em The Universal Gerbe, Dixmier-Douady classes and gauge theory
}. { Lett. Math. Phys.} {\bf 59} 47-60 2002



\bibitem[CMW]{CMW} A.L. Carey, M.K. Murray, B.L. Wang
 {\em Higher bundle gerbes and cohomology classes in gauge theories}.
 J.Geom.Phys. 21 (1997) 183-197.


\bibitem[EM]{EM}
C. Ekstrand, J. Mickelsson
{\it Gravitational anomalies, gerbes, and hamiltonian quantization}
  
  \bibitem[GR]{GR} K.  Gawedski, N.  Reis,{\em  WZW branes and gerbes}.  
  Rev. Math. Phys.  14  (2002),  no. 12, 1281--1334.
  
\bibitem[Getz]{Getz} E. Getzler, {\em The odd Chern character in cyclic homology 
and spectral flow.}    Topology 32 (1993), 489-507.

\bibitem[Hit]{Hit}
N. Hitchin
{\it Lectures on Special Lagrangian Submanifolds}
Winter School on Mirror Symmetry, Vector Bundles
 and Lagrangian Submanifolds (Cambridge, MA, 1999), 
 151--182, AMS/IP Stud. Adv. Math., 23, Amer. Math. Soc., Providence, RI, 2001.

  \bibitem[HS]{HopSin}
   M.J. Hopkins and I.M. Singer,
   {\sl Quadratic functions in geometry, topology,and M-theory}, preprint,
   math.AT/0211216


\bibitem[L]{L} John Lott
{\it Higher-Degree Analogs of the Determinant Line Bundle},  
Comm. Math. Phys.  230  (2002),  no. 1, 41--69.

\bibitem[Loya]{Loya} P. Loya, {\em The Index of $b$-pseudodifferential Operators
 on Manifolds with Corners. } Preprint.

\bibitem[LM]{LM} P. Loya, R. Melrose, {\em Fredholm Perturbations of Dirac Operators on 
Manifolds with Corners.} Preprint.
\bibitem[MaM]{MaM} Private communication.
\bibitem[MaP]{MaP} Mazzeo, R.; Piazza, P. {\em
Dirac operators, heat kernels and microlocal analysis. 
II. Analytic surgery,}  Rend. Mat. Appl. (7)  18  (1998),  no. 2, 221--288.

\bibitem[MP1]{MP1} Melrose, R.B. and Piazza, P.{\em
 Families of Dirac operators, boundaries and the $b$-calculus.}
  J. Diff. Geom.  146, 99--180 (1997) 
\bibitem[MP2]{MP2} Melrose, R.B. and Piazza, P.{\em  An index theorem for families 
of Dirac operators on odd-dimensional manifolds with boundary.}
 J. Diff. Geom. 146, 287--334 (1997) 

\bibitem[Mic]{Mic} J. Mickelsson, {\rm Gerbes, twisted K-theory, and the
supersymmetric WZW model}, hep-th/0206139.

 \bibitem[MR]{MR} J. Mickelsson and S.G. Rajeev.
{\em  Current algebras in d+1
dimensions and determinant line bundles over infinite-dimensional
Grassmanians.},
Commun. Math. Phys. {\bf 116}, 400, (1988).

 \bibitem[Mur]{Mur}
M. K. Murray,
Bundle gerbes,
J. London Math. Soc. (2) {\bf 54}
(1996), no.~2, 403--416.


  \bibitem[MS]{MurSte}
    M. K. Murray, D. Stevenson
  {\sl Higgs fields, bundle gerbes and string structures},
   Comm. Math. Phys. {\bf 243} (2003), no.~3,
   541--555.
 Comm. Math. Phys. {\bf 243} (2003), no.~3,
541--555.


\bibitem[Pi1]{Pi1} P. Piazza, {\em Determinant bundles, manifolds with boundary 
and surgery.}   Comm. Math. Phys.  178  (1996),  no. 3, 597--626.
\bibitem[Pi2]{Pi2} P. Piazza, {\em  Determinant bundles, manifolds with boundary and 
surgery.  II. Spectral sections and surgery rules for anomalies.} 
 Comm. Math. Phys.  193  (1998),  no. 1, 105--124.

\bibitem[PS]{PS}
A. Pressley and G. Segal,
{\it Loop groups}.
Oxford, Clarendon Press, 1986.

\bibitem[Qu]{Qu} D. Quillen, {\em 
Determinants of Cauchy-Riemann operators on Riemann surfaces.}, Func. Anal. Appl. 14,
31-34, (1985).

\bibitem[Sco]{Sco} S. Scott,  {\em Splitting the Curvature of the Determinant Line
Bundle}, Proc. Am. Math. Soc.  128  (2000),  no. 9, 2763--2775.
 
 \bibitem[SW]{SW} S. Scott, K. Wojciechowski, {\em  The $\zeta$-determinant 
 and Quillen determinant for a Dirac operator on a manifold with boundary}, 
  Geom. Funct. Anal.  10  (2000),  no. 5, 1202--1236.
\bibitem[Seg]{Seg} G.B. Segal, {\em Faddeev's anomaly in Gauss' law}
unpublished manuscript, 1986. 
\bibitem[Woj]{Woj} K. Wojciechowski, {\em  The $\zeta$-determinant and 
 the additivity of the $\eta$-invariant on the smooth, self-adjoint Grassmannian}, 
  Comm. Math. Phys.  201, 423--444.
 
\end{thebibliography}
\end{document}